\documentclass[11pt]{article}

\usepackage{latexsym}
\usepackage{amssymb}
\usepackage{euscript}
\usepackage[ansinew]{inputenc}
\usepackage{graphicx}
\usepackage{pstricks}

\let\ssection=\section
\renewcommand{\section}{\setcounter{equation}{0}\ssection}

\def\r{{\mathbb R}}
\def\e{{\mathbb E}}
\def\p{{\mathbb P}}

\newtheorem{theorem}{Theorem}[section]
\newtheorem{proposition}[theorem]{Proposition}
\newtheorem{lemma}[theorem]{Lemma}
\newtheorem{corollary}[theorem]{Corollary}

\def\df{\, {\buildrel{\rm def} \over =}\, }
\def\law{\, {\buildrel{\rm law} \over =}\, }
\def\wcv{\,{\buildrel{(d)} \over\longrightarrow}\,}

\def\was{\,{\buildrel{\rm a.s.} \over\longrightarrow}\, }
\def\weak{\,{\buildrel{w}\over\longrightarrow}\,}
\def\no{\noindent}

\def\as{\rm a.s.}
\def\qed{\hfill$\Box$}
\def\aaa{{\cal A}}
\def\fff{{\cal F}}

\def\norm#1{{\left\| {#1} \right\|}}
\def\arch#1{{\left\{ {#1} \right\}}}

\begin{document}

\vglue30pt

\centerline{\bf\Large  Strong approximations of three-dimensional
Wiener sausages}

\bigskip

\vskip10pt
\vskip10pt

\centerline{Endre C{\small S\'AKI}\footnote{A. R\'enyi Institute
of Mathematics, Hungarian Academy of Sciences, Re\'altanoda u.
13--15, P.O.B. 127, Budapest, H--1364, Hungary. E-mail:
csaki@renyi.hu. Research supported by the Hungarian National
Foundation for Scientific Research, Grants T 037886 and T 043037}
and Yueyun H{\small U}\footnote{Laboratoire de Probabilit\'es et  
Mod\`eles Al\'eatoires (CNRS UMR--7599), Universit\' e
Paris VI, 4 Place Jussieu, F--75252 Paris cedex 05, France.
E-mail: hu@proba.jussieu.fr}}

\bigskip
\bigskip

{\leftskip=2truecm \rightskip=2truecm \baselineskip=15pt

{\small

 {\noindent\slshape\bfseries Abstract.} In this paper we prove
 that the centered three-dimensional Wiener sausage can be strongly
 approximated by a one-dimensional Brownian motion running at a
 suitable  time clock. The strong approximation  gives all possible laws of
 iterated logarithm as well as the convergence in law in terms of process
 for the normalized Wiener sausage. The proof relies on Le Gall
 \cite{LG86b}'s fine $L^2$-norm estimates between the Wiener sausage and the
 Brownian intersection local times.

 \bigskip

 \bigskip

 {\noindent\slshape\bfseries AMS Classification 2000.} 60F15.

} % end of "small"
} % end of "narrower"

\bigskip
\bigskip

 % \newpage

\bigskip
\bigskip
%\tableofcontents
\section{Introduction}

Let $\{B_t, t\ge0\}$ be a $d$-dimensional Brownian motion and let
$K \subset \r^d$ be a compact set. The Wiener sausages associated
to $(B_t)$ and $K$ are the compact sets $$ S^K(0, t) \df
\bigcup_{0\le s\le t} (B_s + K), \qquad t\ge0.$$

Denote by $m(dx)$ the Lebesgue measure. The volume of Wiener
sausage  $m(S^K(0,t))$  is a very rich topic of researches, we
only mention some recent references (\cite{CL04}, \cite{H02},
\cite{HK01}), see  their lists of references for related studies.

We restrict our attention to the three dimensional case ($d=3$).
Denote by ${\cal C}(K)$ the Newtonian electrostatic capacity of
the compact $K$, see \cite{IK74} pp. 250.  Concerning the
asymptotic of $m(S^K(0, t))$ when $t \to \infty$, Kesten, Spitzer
and Whitman (\cite{IK74} pp. 252, Problem 4), Spitzer \cite{S64},
and Le Gall \cite{LG88} proved respectively the following
results (\ref{1}), (\ref{2}) and ({\ref{3}):

 \noindent {\bf Theorem A} {\it
When $t\to\infty$,
\begin{eqnarray}
{1\over t} m(S^K(0, t))  & \was &   {\cal C}(K),  \label{1} \\
  \e \Big(m(S^K(0, t)) \Big) & =& {\cal C}(K) t + {4\over (2\pi)^{3/2}}
 {\cal C}(K)^2 \, \sqrt {t} +  o(\sqrt t), \label{2} \\
 {m(S^K(0, t)) - {\cal C}(K) t \over \sqrt{ t\log t}} & \wcv & { {\cal C}(K)^2
 \over \pi \sqrt  2 }  \, {\cal N}, \label{3}
 \end{eqnarray}
where $\cal N$ denotes a centered  real-valued Gaussian variable,
of variance $1$, and ${\cal C}(K)$ denotes   Newtonian  capacity.
}

This paper aims at establishing  a strong approximation of the
process $\{m(S^K(0,t)), t\ge0\}$:

\begin{theorem}\label{T1}  Let $l\ge1$ and $K_1, ..., K_l$ be $l$ compact sets in $\r^3$.  There exists a  one dimensional Brownian
motion $(\beta(t), t\ge0)$ such that almost surely, for all $1\le
i \le l$,
$$ m(S^{K_i}(0,t))= {\cal C}(K_i) t + { {\cal C}(K_i)^2 \over \pi \sqrt 2} \,
\beta(t \log t) + o \big( t^{1/2}\, (\log t)^{1/4 + \delta}\big),
$$ for any positive constant $\delta>0$.
\end{theorem}

The above strong approximation is in agreement with Bass and
Kumagai \cite{BK00} where they deal with the range of  a
three-dimensional (and two-dimensional) random walk. Note also
that the Brownian motion $\beta$ can be   chosen simultaneously
the same  for all the $l$ compact sets $(K_i)_{1\le i\le l}$, a
fact   already   pointed out by Le Gall \cite{LG88}.

Theorem \ref{T1} implies the convergence in terms of process of
the normalized Wiener sausage to Brownian motion. The error term $
o \big( t^{1/2}\, (\log t)^{1/4 + \delta}\big)$ being sufficiently
small, in addition to weak convergence, we can deduce functional almost
sure central limit theorem and laws of the iterated logarithm.

\begin{corollary} Let $K\subset \r^3$ be a compact set such that
${\cal C}(K)>0$. Put
$$
M_t(x)={\pi \sqrt{2}\over {\cal C}(K)^2\sqrt{ t \log t}} \, \left(
m(S^{K}(0,xt))-{\cal C}(K) xt\right),\quad 0\le x\le 1.
$$
We have
\begin{eqnarray*}
{\cal L}\left(M_t(\cdot)\right)&\weak& {\cal L}(\beta),\quad t\to\infty\\
{1\over \log T}\int_1^T{1\over t}\delta\left(M_t(\cdot)\right)\, dt
&\weak& {\cal L}(\beta), \quad{\rm a.s.},\quad T\to\infty\\
\limsup_{t\to\infty}
{1\over \sqrt{ t \log t \log\log t}} \left(m(S^{K}(0,t))-{\cal
C}(K) t\right) &= & { {\cal C}(K)^2 \over \pi  }, \qquad \as \\
\liminf_{t\to\infty}  \left(  \log \log t\over t \log t
\right)^{1/2}\, \sup_{0\le s\le t} \big| m(S^{K}(0,s))-{\cal C}(K)
s\big| &=& { {\cal C}(K)^2 \over 4  }, \qquad \as,
\end{eqnarray*}
where $\weak$ denotes weak convergence in $C[0,1]$, ${\cal L}(\cdot)$
denotes the law of the process in bracket and $\delta(f(x))$
stands for a point mass at $f\in C[0,1]$.
\end{corollary}

The proof of Theorem \ref{T1} relies on two steps:

(i)  Establish the strong approximation of  the Wiener sausage  by
an intersection local time;

(ii) Use a Tanaka-Rosen-Yor (\cite{R86}, \cite{Y85}) type formula
and  approximate the intersection local time by a Brownian motion.

This paper is organized as follows: In Section 2, we give a
general estimate on  the rate of growth of a  process through its
$L^p$-norms; Section 3 is devoted to the intersection local times,
in particular to show the above step (i); the step (ii) is done in
Section 4. Finally, we prove Theorem 1.1 in Section 5.

Throughout this paper,   $\norm{\xi}_p$ denotes the $L^p$-norm of
a random variable $\xi$, while  $c_p, c'_p, c''_p$ denote some positive
constants, depending on $p>0$, whose exact values may vary from
one paragraph to another.

\section{Rate of growth}

We shall estimate  the rate of growth of some process through its
$L^p$-norms. The following result  on the modulus of continuity,
due to Barlow and Yor \cite{BY82}, is obtained as a consequence
of Garsia-Rodemich-Rumsey lemma:

\begin{lemma}\label{L:barlowyor} {\bf (\cite{BY82}, (3.b))} Let $r>0$. Assume that  $(\xi_t,
0< t < r)$ is a real-valued  continuous process such that for some
constants $a>0,  b>1$ and $\kappa>0$, $$ \e  \left| \xi_t - \xi_s
\right|^a \, \le\, \kappa \, |t-s|^b, \quad \forall\, 0<s\le
t<r.$$ Then for any $0< \gamma < b -1$, there exists some constant
$c(a, b, \gamma)>0$
 such that $$ \e \left( \sup_{0< s \not= t<r}{ |\xi_t - \xi_s
|^a \over |t-s|^{\gamma}}\right) \le c(a,b, \gamma) \, \kappa\,
r^{b-\gamma}.$$
\end{lemma}

\begin{proposition}\label{P:inc1} Let $(X_t, t>0)$ be a
real-valued  process. Let $0< a \le1$, $b_2\ge b_1 \ge0$ be three
constants. Assume that for any $p>1$, there exists some constant
$c_p>0$ such that
\begin{eqnarray} \left\|  X_t  \right\|_p  &\le& c_p  \, t^a\, (\log
t)^{b_1}, \quad \forall \, t>t_0,\label{pi1} \\
  \norm{  X_t  - X_s }_p  &\le &c_p  \, (t-s )^a\, (\log
t)^{b_2},\quad   t_0\le s\le t-1, \label{pi2}\\
   \norm{ \sup_{t\le u\le t+1} |X_u - X_t|}_p & \le & c_p \, t^{a/2}, \quad
   \forall\, \, t>t_0,  \label{pi3} \end{eqnarray}
   for some constant $t_0>1$. Then for any $\epsilon>0$, $$ X_t = o\left( t^a
(\log t)^{b_1+ \epsilon}\right),  \qquad {\rm a.s.}, \, \, t\to\infty.$$
\end{proposition}

The power $a/2$ in (\ref{pi3}) can be replaced by any positive
constant smaller than $a$. Note that in (\ref{pi2}), we only
require for the increments $X_t - X_s$  when $t-s >1$, the local
fluctuations are controlled by (\ref{pi3}).

{\no\bf Proof:} Let $p> 4/a$ and $n> t_0$.  By (\ref{pi3}),
\begin{eqnarray*} \p\Big( \max_{t_0\le j \le n} \sup_{j \le t\le
j+1} |X_t - X_j| > n^a\Big) &\le & \sum_{j=[t_0]}^n \p\Big(  \sup_{j
\le t\le j+1} |X_t - X_j|
>n^a\Big)
     \\&\le & \sum_{j=[t_0]}^n  n^{-ap}\, c_p^p\, (1+n)^{ap/2} \\
     &\le & c'\, n^{- ap/2 +1},
\end{eqnarray*}

\no whose sum on $n$ converges. Applying  Borel-Cantelli's lemma
shows that $$ \max_{0\le j \le n} \sup_{j \le t\le j+1} |X_t -
X_j|  = O(n^a), \qquad {\rm a.s.}, \, \, n\to\infty.$$

\no Then it remains to prove that$$ X_n= o(n^a\, (\log
n)^{b_1+\epsilon}), \qquad {\rm a.s.}, \,\, n\to\infty.
$$

 Let $0< \epsilon < a/4$ and fix any constant $\eta \in
(0,1/2) $  satisfying  $\eta b_2 \le a( 1- \eta)/2$.  Let $p>
8/(\eta \epsilon)$. Consider large $j$ and define $n_j=
[e^{j^\eta}]$. Since $\epsilon>0$ is arbitrary, it suffices to
prove that
\begin{equation} \label{ic36}\sum_j \p\Big( \max_{n_j \le k\le
n_{j+1}} |X_k| \ge n_j^a \, (\log n_j)^{b_1+ \epsilon}\Big) <
\infty. \end{equation}

To this end, we have by triangle inequality
\begin{eqnarray*} && \p\Big( \max_{n_j \le k\le n_{j+1}}
|X_k| \ge n_j^a \, (\log n_j)^{b_1+ \epsilon}\Big) \\
     &\le & \p\Big( |X_{n_j} |  \ge  {1\over2} n_j^a \, (\log n_j)^{b_1+
\epsilon}\Big) + \p\Big( \max_{n_j \le k\le n_{j+1}} |X_k -
X_{n_j} | \ge {1\over2} n_j^a \, (\log n_j)^{b_1+
\epsilon}\Big) \\
    & \df  & I_1+ I_2.
\end{eqnarray*}

 Using the $L^p$-norm of
$X_{n_j}$,
\begin{eqnarray*} I_1 &\le & \left({1\over2} n_j^a \, (\log
n_j)^{b_1+ \epsilon}\right)^{-p} \, c_p^p\,n_j^{ap}\, (\log
n_j)^{p b_1}
    \\ &\le &  c'_p\, j^{-  p\epsilon  \eta } \le j^{-8}.
    \end{eqnarray*}

We shall  apply Lemma \ref{L:barlowyor} to estimate $I_2$.
Consider a continuous process $\xi$ such that $ \xi_k= X_{k+ n_j}
- X_{n_j}$ for $0\le k \le n_{j+1} - n_j$, and $\xi_\cdot$ is
linear on each interval $[k, k+1]$. From (\ref{pi2}), we easily
deduce that for any $0<u, v < n_{j+1} - n_j$,  $$ \e  |\xi_u - \xi_v|^p  \le  \left\{%
\begin{array}{ll}
     c_p^p \, |v-u|^{ap}\, (\log  n_{j+1})^{ p
b_2}, & \hbox{ $|u-v| >2$;} \\
      \\
    c'_p \, |v-u|^p\, (\log  n_{j+1})^{ p
b_2} , & \hbox{$|u-v| <1.$} \\
\end{array}%
\right.
$$

Since $a\le1$, there exists  some constant $  c_{a,p}>0$ such that
  for
any $0< u, v <n_{j+1} - n_j$, $$ \e | \xi_u - \xi_v |^p \le \,
c_{a,p}\, (\log  n_{j+1})^{b_2 p}\, |u-v|^{ap}.$$

Applying Lemma \ref{L:barlowyor} with $\gamma = p(a - \epsilon/2)
< a p-1$ gives that for some constant $c(\epsilon, p)>0$
$$ \e\Big( \sup_{ 0\le u \not= v<n_{j+1}-n_j} { |\xi_v -\xi_u|^p\over |v-u|^{p(a - \epsilon/2)}} \Big) \,\le \,
c(\epsilon, p)\, (n_{j+1}- n_j)^{p\epsilon/2} (\log  n_{j+1}
)^{b_2 p}.$$

It follows from Chebyshev's inequality that \begin{eqnarray*} I_2
&\le& \p\Big( \sup_{ 0< k \le  n_{j+1} - n_j} { |\xi_k|^p\over
k^{p(a - \epsilon/2)}}
> 2^{-p}\, { n_j^{ap}\over (n_{j+1}- n_j)^{p (a-\epsilon/2)}}
(\log n_j)^{p(b_1+\epsilon)}\Big) \\
    &\le & c_{  \epsilon, p} 2^p \, \left({n_{j+1}- n_j\over
    n_j}\right)^{ap}  \, (\log n_{j+1})^{p b_2}\, (\log
    n_j)^{-p(b_1+\epsilon)} \\
    &\le &  c'\, j^{ - (1-\eta) p a  -
    p\eta (b_1+\epsilon) + p \eta b_2} \\
    &\le & c'\, j^{- a p(1-\eta)/2 } \le  j^{-2},
\end{eqnarray*}
by our choice of $p$ and $\eta$. Then the sum of $I_1$ and $I_2$
over $j$ converges and it proves (\ref{ic36}) hence the Proposition.
\hfill$\Box$

When the estimate (\ref{pi2}) on the increments holds for all
$t>s$, we may get rid of the condition on the local fluctuations
(\ref{pi3}), as stated in the following result:

\begin{proposition}\label{P:inc2} Let $(X_t, t>0)$ be a
real-valued  process. Let $0< a \le1$, $b_2\ge b_1 \ge0$ be three
constants. Assume that for any $p>1$, there exists some constant
$c_p>0$ such that
\begin{eqnarray*} \left\|  X_t  \right\|_p  &\le& c_p  \, t^a\, (\log
t)^{b_1}, \quad \forall \, t\ge t_0, \\
 \left\|  X_t  - X_s \right\|_p  &\le & c_p  \, (t-s )^a\, (\log
t)^{b_2},\quad   t_0\le s\le t,
\end{eqnarray*} for some $t_0>1$. Then for any $\epsilon>0$, $$ X_t = o\left( t^a
(\log t)^{b_1+ \epsilon}\right), \qquad {\rm a.s.},\, \, t\to\infty.$$
\end{proposition}

{\no\bf Proof:} The proof goes in the same way as in Proposition
\ref{P:inc1}. Let us only mention the main difference. Let $
\epsilon , \eta  , p , n_j$ be the same as in Proposition
\ref{P:inc1}. To control the probability
$$ \p\Big( \sup_{n_j \le
t\le n_{j+1}} |X_t - X_{n_{j+1}} | \ge {1\over2} n_j^a \, (\log
n_j)^{b_1+ \epsilon}\Big),$$
we define  $\xi_u := X_{ u
(n_{j+1}-t_0) + t_0}, 0\le u\le 1 $. Then for any $0< u, v <1$, $$
\e | \xi_u - \xi_v |^p \le \, c_p^p\, n_{j+1}^{ap} \, (\log
 n_{j+1})^{b_2 p}\, |u-v|^{ap}.$$

Applying Lemma \ref{L:barlowyor} with $\gamma = p(a - \epsilon/2)
< a p-1$ gives that for some constant $c(\epsilon, p)>0$
$$ \e\Big( \sup_{ 0< s \not= t<n_{j+1}} { |X_t -X_s|^p\over |t-s|^{p(a - \epsilon/2)}} \Big) \,\le \,
c(\epsilon, p)\, n_{j+1}^{p\epsilon/2} (\log  n_{j+1} )^{b_2 p}.$$
Thus    the proposition follows in the same way. \hfill$\Box$

\section{Intersection local times}

Denote by $\{\alpha(x, F), x\in \r^3, F  \mbox{ bounded measurable
}   F \subset \r_+^2\}$ the family of intersection local times of
$B$. For any $F$, almost surely $\alpha(x, F)$ is the density of
the measure $  \int_F d s du \delta_{(B_s- B_u)}(dx)$, which means
that for any test function $f: \r^3 \to \r_+$,
\begin{equation} \label{ilt} \int_{(s,u) \in F} \, d s d u\, f(B_s
-B_u) = \int_{\r^3} d x f(x) \alpha(x, F). \end{equation}

We refer to Rosen \cite{R83}, Le Gall \cite{LG85} for the
construction of the family $\{\alpha(x, F)\}$. We shall be
particularly interested in two families of subsets $F$:
\begin{eqnarray*} \aaa_t &=& \big\{(s, u): 0\le s\le  u-1 \le
t-1\big\}, \qquad t>1.\\
       \fff_\delta &=& \big\{ (s,u): 0\le s \le u - \delta \le
       1-\delta\big\}, \qquad 0<\delta <1.
\end{eqnarray*}

Let us introduce the following notation: For $0\le a < b$,
$S^K_\epsilon(a,b)\df \bigcup_{a\le t\le b} ( B_t+ \epsilon K)$
and $ S^K_\epsilon = \bigcup_{0\le s \le 1} (B_s + \epsilon K)$
when there is no risk of confusion.  Then we have the following
scaling property: For $t>1$, $$ \Big( m(S^K(0,t)), \alpha(0,
\aaa_t)\Big) \law \Big( t^{3/2} m(S^K_{1/\sqrt t}), t^{1/2}
\alpha(0, \fff_{1/t})\Big).$$

 Write $$ \arch{\xi} = \xi - \e (\xi) .$$ The
following $L^p$-norm estimate is essentially due to Le Gall who
proved the case $p=2$ in \cite{LG88}, Theorem 3.1. His arguments
can be adopted to deal with the $L^p$ case.

\begin{proposition}\label{P:lp}  Fix a compact set $K\subset \r^3$. 
For any $p\ge2$, there exists some constant $c_p>0$ depending on $p$ and 
$K$ such that for all $0<\epsilon < 1/2$,
$$ \norm{ {1\over \epsilon^2} \arch{  m(S^K_\epsilon) } - {\cal C}(K)^2
\arch{ \alpha(0, \fff_{\epsilon^2})  }  }_p \le c_p .$$
Consequently, for any $t>1$,
$$ \norm{   \arch{ m(S^K(0,t)) } - {\cal C}(K)^2
\arch{ \alpha(0, \aaa_t)  }  }_p \le c_p  \, \sqrt t.$$
\end{proposition}

We shall make use of the following Rosenthal inequality, see
Petrov \cite{P95}, Theorem 2.9:

\begin{lemma}\label{L:ros}   For any  $n$   centered, independent and
identically distributed  random variables $X_1,...X_n$
there exists a constant $c_p>0$ depending only on
$p\ge2$ such that
\begin{equation}\norm{  \sum_{i=1}^n X_i }_p  \le  \left\{%
\begin{array}{ll}c_p \,
\sqrt{n}\,
\norm{ X_1}_p ,  \\
     \\
    c_p \, \Big( \sqrt{n}\, \norm{
X_1}_2 + n^{1/p}\, \norm{  X_1 }_p\Big).   \\
\end{array}%
\right.  \label{ros} \end{equation}
\end{lemma}

{\no\bf Proof of Proposition \ref{P:lp}:} The proof of the case
$p>2$ is almost the same as the case $p=2$ given by Le Gall
(\cite{LG88}, pp. 1006). In fact, let $n$ be an integer such that
$1/4 < \epsilon 2^{n/2} \le 1/2$. Then using the   notation
$S^K_\epsilon(a,b)\df \bigcup_{a\le t\le b} ( B_t+ \epsilon K)$,
$I_-^{(k,q)}=({2q-2\over 2^k}, {2q-1\over 2^k})$ and
$I_+^{(k,q)}=({2q-1\over 2^k}, {2q\over 2^k})$ for $1\le q\le
2^{k-1}$,
$$ \arch{ m(S^K_\epsilon) } = \sum_{k=1}^{2^n}\, \arch{ m\left(
S^K_\epsilon\left({k-1\over 2^n}, {k\over 2^n}\right)\right) } -
\sum_{k=1}^{n}
\sum_{q=1}^{2^{k-1}} \, \arch{ m(S^K_\epsilon(I_-^{(k,q)}) \cap
 S^K_\epsilon(I_+^{(k,q)}) )}.$$

Therefore by the independent increment property of Brownian
motion, we can apply (\ref{ros}) and get
\begin{eqnarray}
\norm{  {1\over \epsilon^2}  \sum_{k=1}^{2^n}\, \arch{ m\left(
S^K_\epsilon\left({k-1\over 2^n}, {k\over 2^n}\right)\right) } }_p &\le &
c_p\epsilon^{-2} \, 2^{n/2}\, \norm{ \arch{ m\left(S^K_\epsilon\left(0,
{1\over 2^n}\right)\right)  } }_p \nonumber
    \\ &\le & c_p \epsilon^{-2} 2^n \sup_{1/4 < r\le 1/2}
    \norm{ \arch{ S^K_r(0,1)} }_p  \nonumber\\
        &\le & c'_p, \label{12}
\end{eqnarray}

\no where the second inequality follows from the  scaling
property:  $m(S^K_\epsilon(0, {1\over 2^n})) \law 2^{-3n/2} \,
m(S^K_r(0,1))$ with $r= \epsilon 2^{-n/2} \in ({1\over4},
{1\over2}]$. On the other hand,
\begin{eqnarray*}&& \norm{
{1\over \epsilon^2}\sum_{k=1}^{n} \sum_{q=1}^{2^{k-1}} \, \arch{
m(S^K_\epsilon(I_-^{(k,q)}) \cap
 S^K_\epsilon(I_+^{(k,q)}) )}- {\cal C}(K)^2 \arch{
\alpha(0, {\cal F}_{2^{-n}})} }_p \\
    &\le& \sum_{k=1}^n \, \norm{
{1\over \epsilon^2} \sum_{q=1}^{2^{k-1}} \, \arch{
m(S^K_\epsilon(I_-^{(k,q)}) \cap
 S^K_\epsilon(I_+^{(k,q)}) )} - {\cal C}(K)^2 \arch{
\alpha(I_-^{(k,q)} \times I_+^{(k,q)}) } }_p \\
    &\le& c_p\, \sum_{k=1}^n \, \Big( 2^{(k-1)/2} \norm{
    \xi(\epsilon, k) }_2+ 2^{(k-1)/p} \norm{
    \xi(\epsilon, k) }_p \Big),
\end{eqnarray*}

\no by applying the second inequality in (\ref{ros}) and where $$
\xi(\epsilon, k) = {1\over \epsilon^2}\, \arch{ m(S^K_\epsilon(0,
2^{-k}) \cap
  S^K_\epsilon(2^{-k}, 2^{-k+1})) } - {\cal C}(K)^2 
\arch{  \alpha((0, 2^{-k})\times (2^{-k}, 2^{-k+1}))}$$

Le Gall (\cite{LG88}, pp.1006) has already shown that
\begin{equation} \label{14} \sum_{k=1}^n \, 2^{(k-1)/2} \norm{
\xi(\epsilon, k) }_2 = O(1).
\end{equation}

We shall prove that $$ \sum_{k=1}^n \, 2^{(k-1)/p} \norm{
    \xi(\epsilon, k) }_p = O(1),$$ which together with (\ref{12})
    and (\ref{14}) implies the Proposition.

    To this end, by scaling, \begin{eqnarray*}
 \sum_{k=1}^n \, 2^{(k-1)/p} \norm{
    \xi(\epsilon, k) }_p &=& 2^{-1/p} \, \sum_{k=1}^n 2^{ - (
    {1\over 2} - {1\over p}) k}\, 
    \norm{  r_k^{-2} \arch{ m(S_{r_k}^K \cap \widetilde S_{r_k}^K ) } - 
    {\cal C}(K)^2 \arch{{\cal L}_2([0,1]^2)} }_p
    \\ &\le & 2^{-1/p+1} \, \sum_{k=1}^n 2^{ - (
    {1\over 2} - {1\over p}) k}\, \left(\norm{  r_k^{-2}  m(S_{r_k}^K 
    \cap \widetilde S_{r_k}^K)}_p+ {\cal C}(K)^2\, 
    \norm{   {\cal L}_2([0,1]^2)}_p\right),
    \end{eqnarray*}
where $r_k= \epsilon 2^{k/2}<1$, $S_{r_k}^K= S_{r_k}^K(0,1)$ and
${\cal L}_2([0,1]^2)$ denotes the intersection local times of two
independent Brownian motion $B$ and $\widetilde B$ over $[0,1]^2$.
Using Le Gall (\cite{LG86b}, Corollary 3.2), we have that for any
compact $K$, $\norm{{\cal L}_2([0,1]^2) }_p<\infty$ and
$$ \sup_{0<r<1}\norm{  r^{-2}   m(S_r^K \cap \widetilde
S_r^K )}_p <\infty.$$

\no It follows that  $$ \sum_{k=1}^n \, 2^{(k-1)/p} \norm{
    \xi(\epsilon, k) }_p \le  c_p\,\sum_{k=1}^n 2^{ - (
    {1\over 2} - {1\over p}) k} =O(1) ,  $$ completing the proof.  \hfill$\Box$

The main technical estimate in this section is the following
$L^p$-norm:

\begin{lemma}\label{L:m2b} Let $p\ge 1 $. Denote by $D$ the unit ball of $\r^3$.
We have for all $t>1$, $$ || m\big( S^D(0,t) \cap \widetilde
S^D(0, \infty)\big) ||_p \le c_p \, \sqrt{t},$$ where $c_p>0$ only
depending on $p$ and $\widetilde S^D$ denotes the Wiener sausage
associated with another Brownian motion independent of $B$.
\end{lemma}

We may choose $c_p= c^p$ in the above lemma  for some numerical
constant $c>1$.

{\no\bf Proof of Lemma \ref{L:m2b}:} By scaling $$  m\Big(
S^D(0,t) \cap \widetilde S^D(0, \infty)\Big) \,\law  \,
\epsilon^{-3}\,  m\Big( S^D_\epsilon(0,1) \cap \widetilde
S^D_\epsilon(0, \infty)\Big),$$ with $\epsilon = t^{-1/2}$. It is
equivalent  to prove that for any integer $p\ge 1$,
\begin{equation} \e\Big( m\Big( S^D_\epsilon(0,1) \cap \widetilde
S^D_\epsilon(0, \infty)\Big) \Big)^p  \le c'_p  \, \epsilon^{2p},
\qquad 0< \epsilon<1. \label{eps2p}\end{equation}

Denote by $T^K_\epsilon(x)$ the first entry time into $x- \epsilon
K$: $$ T^K_\epsilon(x) = \inf\{t\ge0: B_t \in x - \epsilon K\},$$
with convention that $\inf \emptyset =\infty$.  Define $$ \psi(r)
= {1\over r} 1_{(r<1)} + e^{ - r^2/16}, \qquad r>0.$$ Recall the
following estimate, see Le Gall \cite{LG88}, Lemma 3.2:
\begin{eqnarray} \p\Big( T_\epsilon^D(x)  \le 1\Big) & \le & c\,
\psi( |x|), \label{tdx1} \\
\p\Big(T_\epsilon^D(x)  < \infty \Big) & =  &  
{ \epsilon \over  |x|} \wedge 1.  \label{tdx2}
\end{eqnarray}

\no Note that \begin{eqnarray}\label{eps3p} \Upsilon_p(\epsilon)
  &\df& \e \Big( m\big( S^D_\epsilon(0,1) \cap \widetilde
S^D_\epsilon(0, \infty)\big) \Big)^p  \nonumber \\
    &= & \int d
x_1...dx_p \, \p\Big( \max_{1\le i\le p} T_\epsilon^D(x_i) \le1
\Big) \, \p\Big( \max_{1\le i\le p} T_\epsilon^D(x_i) <
\infty\Big) \nonumber \\
    &=& \int_{ \exists i\not= j: |x_i- x_j| \le 2 \epsilon} +  
\int_{ \forall i\not= j: |x_i- x_j| > 2 \epsilon}, \end{eqnarray}

\no with respect to the obvious density.  By symmetry, we have
\begin{eqnarray}  \int_{ \exists i\not= j: |x_i- x_j| \le 2 \epsilon} &\le&
  {p(p-1)\over2} \,  \int d x_1... d x_{p-1} \, 
  \int_{ |x_p - x_{p-1}|\le 2\epsilon} d x_p  \nonumber\\
  &\le & c \,(2\epsilon)^3\, {p(p-1)\over2} \, \int d x_1... d x_{p-1} \,
 \p\Big( \max_{1\le i\le p-1} T_\epsilon^D(x_i)
\le1 \Big) \, \p\Big( \max_{1\le i\le p-1} T_\epsilon^D(x_i) <
\infty\Big)  \nonumber \\
    &\le& 4 c p^2 \epsilon^3\, \Upsilon_{p-1}(\epsilon).
    \label{u11}
\end{eqnarray}

Let ${\cal S}_p$ be the set of permutations on $[1,...,p]$.  By
the strong Markov property, we have
\begin{eqnarray*}  &&\p\Big(\max_{1\le i\le p} T_\epsilon^D(x_i)< \infty) \\ &
\le & \sum_{\sigma \in {\cal S}_p} \, \p\Big( T_\epsilon^D(
x_{\sigma(1)}) \le ...\le T_\epsilon^D( x_{\sigma(p)}) <
\infty\Big) \\
    &\le & \sum_{\sigma \in {\cal S}_p} \, \p\Big( T_\epsilon^D(
x_{\sigma(1)}) \le ...\le T_\epsilon^D( x_{\sigma(p-1)}) <
\infty\Big)  \sup_{ |y- x_{\sigma(p-1)}| \le \epsilon} \p_z
\Big(T_\epsilon^D( x_{\sigma(p)}) < \infty\Big)  \\
    &\le & \sum_{\sigma \in {\cal S}_p} \, \prod_{i=1}^p  
    \Big({ \epsilon \over ( |x_{\sigma(i)} -
    x_{\sigma(i-1)}| -\epsilon)^+} \wedge 1\Big),
\end{eqnarray*}

\no by means of (\ref{tdx2}), with convention that $\sigma(0)=0$
and $x_0=0$. Similarly, we get from (\ref{tdx1}) that
$$ \p\Big(\max_{1\le i\le p} T_\epsilon^D(x_i) \le 1 \Big) \le
(c_p)^p\, \epsilon^p\, \sum_{\sigma \in {\cal S}_p}  \prod_{i=1}^p
 \psi( |x_{\sigma(i)} -    x_{\sigma(i-1)}| -\epsilon)^+).$$

\no Plugging these into (\ref{eps3p}), we obtain that  $$\int_{
\forall i\not= j: |x_i- x_j| > 2 \epsilon}  d x_1...dx_p \,
\p\Big( \max_{1\le i\le p} T_\epsilon^D(x_i) \le1 \Big) \, \p\Big(
\max_{1\le i\le p} T_\epsilon^D(x_i) < \infty\Big) \le  c^p\, p!\,
\epsilon^{2p}\, J_p^*,$$ with    \begin{eqnarray*} J_p^* &\df&
\max_ {\sigma \in {\cal S}_p} \,  J_p(\sigma) \df \max_ {\sigma
\in {\cal S}_p} \,  \int dx_1...dx_p\, \prod_{i=1}^p
 {\psi( |x_{\sigma(i)} -    x_{\sigma(i-1)}|  )  \over  |x_i -x_{i-1}| }
 \Big).
\end{eqnarray*}

We are going to prove that there exists some constant $c>0$ such
that
\begin{equation} \label{pint1} J_p^* \le c^p \, 2^{3p^2} ,
\end{equation}

\no  which in view of (\ref{eps3p}) and (\ref{u11}) implies
(\ref{eps2p}) and completes the proof.

  To show
(\ref{pint1}), we firstly remark that the function $y (\in \r^3)
\to \int_{\r^3} { \psi(|x|)   \over |x -y|}  dx $ is continuous on
the whole $\r^3$ and goes to $0$ when $ |y|\to \infty$, hence
\begin{equation} \label{py1}  \sup_{y\in \r^3} \int_{\r^3}  
{ \psi(|x|)  \over |x -y|} d x = c < \infty. \end{equation}

This implies that for any $y, z \in \r^3$, \begin{eqnarray} \int d
x { \psi(|x|) \, \psi(|x-y|) \over |x-z|} & = &  \int_{ |x-y| \le
|y|/2} + \int_{ |x-y| >  |y|/2 }  \nonumber
\\ &\le &  \psi({ |y|\over 2})\int_{ |x-y| \le
|y|/2} { \psi(|x-y|)\over |x-z|} d x +   \psi({ |y|\over 2})\int_{
|x-y| >  |y|/2 } { \psi(|x|)\over |x-z|} d x   \nonumber\\
    &\le & 2 \,c \, \psi({ |y|\over 2}) , \label{py2}
\end{eqnarray}

\no where in the first inequality, we make use of the monotonicity
of $\psi(\cdot)$ and the fact that $|x| \ge |y|/2$ on $\{ |x-y|\le
|y|/2\}$. We prove (\ref{pint1}) by recurrence on $p$. When $p=1$,
$J_1^*= \int {\psi(|x|) \over |x|} d x < \infty$. Let $p\ge2$ and
consider $\sigma \in {\cal S}_p$. If $\sigma(p)=p$, we have
\begin{eqnarray*}
J_p(\sigma) &=& \int   d x_1... d x_{p-1} \prod_{i=1}^{p-1}
 {\psi( |x_{\sigma(i)} -    x_{\sigma(i-1)}|  )  \over  |x_i -x_{i-1}| }
 \Big)\, \int d x_p { \psi( |x_p - x_{\sigma_{p-1}}|) \over |x_p -
 x_{\sigma_{p-1}}|} \\
    &\le & c\, J_{p-1}(\sigma) \le c\, J_{p-1}^*,
\end{eqnarray*}

\no by (\ref{py1}). If $\sigma(p)\not=p$, we denote by
$j=j(\sigma,p) \le p-1$ such that $\sigma(j)=p$. We have
\begin{eqnarray*} J_p(\sigma) &=& \int   
{ d x_1... d x_{p-1} \over |x_1|...|x_{p-1} - x_{p-2}|}\,  
\prod_{i=1, i \not= j-1,  i \not=j}^{p }
 \psi( |x_{\sigma(i)} -    x_{\sigma(i-1)}|  )  \, \int  d x_p\,
 { \psi( |x_{\sigma_j} - x_{\sigma_{j-1}}|) \psi( |x_{\sigma_{j+1}} -
 x_{\sigma_j}|)\over |x_p -x_{p-1}|}   \\
 &\le & 2c \, \int   { d x_1... d x_{p-1} 
 \over |x_1|...|x_{p-1} - x_{p-2}|}\,  \prod_{i=1, i \not= j-1,  i 
 \not=j}^{p }\psi( |x_{\sigma(i)} -    x_{\sigma(i-1)}|  )  \,\,  \psi(  { |
 x_{\sigma_{j+1}}-  x_{\sigma_{j-1}}| \over2} )  \\
    & \le & 2 c \int   d x_1... d x_{p-1} \prod_{i=1}^{p-1}
 {\psi(  |x_{\tilde\sigma(i)} -    x_{\tilde \sigma(i-1)}| /2  )  
\over  |x_i -x_{i-1}| } \,
\end{eqnarray*}

\no where we used (\ref{py2}) for the first inequality and the
monotonicity of $\psi$ for the second inequality,  the new
permutation $\tilde\sigma \in {\cal S}_{p-1}$ is defined by:
$\tilde\sigma(k)= \sigma(k) $ for $k<j$ and $\tilde\sigma(k)=
\sigma(k+1)$ for $j\le k \le p-1$. By the change of variable $x= 2
\tilde x$, we get $$ J_p(\sigma) \le 2c \, 2^{3(p-1)}\,
J_{p-1}(\tilde\sigma) \le 2^{3p-2}\, c\, J^*_{p-1}.$$

\no Hence $J^*_p \le 2^{3p-2}\, c\, J^*_{p-1} \le 2^{3p^2} c^p$,
yielding that for all $p\ge1$ and $\epsilon<1/2$, $$\e\Big( m\Big(
S^D_\epsilon(0,1) \cap \widetilde S^D_\epsilon(0, \infty)\Big)^p
\Big)  \le \epsilon^{2p}\, c \, 2^{4 p^2}.$$ \qed

\begin{corollary}\label{C:inc} Fix the compact set $K$. For any $p>1$, 
there exists some constant $c_p >0$ such that  \begin{eqnarray}    \norm{ 
\arch{m(S^K(0, t))} - \arch{ m(S^K(0, s))}  }_p  & \le  &
    c_p\, \sqrt{(t-s) \log t },  \, \, 0<s<t-1 ;   \label{inc1} \\
  \norm{ \sup_{t\le
u\le t+1} \big| m(S^K(0, u)) -   m(S^K(0, t)) \big|  }_p  & \le
 & c_p. \label{inc2}
     \end{eqnarray}
\end{corollary}

{\no\bf Proof:} Consider firstly the case $t-s>1$. Write
$F\backslash G:= F\cap G^c$. Define $\widetilde B_u:= B_{u+s} -
B_s, u\ge0$
 and $\widetilde S^K$ the associated Wiener sausage. We have 
\begin{eqnarray*} \arch{ m(S^K(0, t))} - \arch{
m(S^K(0, s))}  &=& \arch{ m(S^K(s,t) \backslash  S^K(0, s))} \\
    & \law & \arch{ m( \widetilde S^K(0,t-s) \backslash  S^K(0, s))} \\
    &=& \arch{ m( \widetilde S^K(0,t-s))} - \arch{ m( \widetilde
    S^K(0,t-s) \cap S^K(0, s))}.
\end{eqnarray*}

Firstly, the convergence in law of (\ref{3}) holds in fact in
$L^p$ for any $p>1$. This together with (\ref{2})  shows that 
$$
\norm{ \arch{ m( \widetilde S^K(0,t-s))} }_p \le c_p\, \sqrt{
(t-s) \log (3+t-s) }, \qquad \forall\, 0<s<t-1.$$

Secondly, we have  \begin{eqnarray*} \norm{
 \arch{ m( \widetilde     S^K(0,t-s) \cap S^K(0, s))} }_p &\le & 2
 \norm{    m( \widetilde     S^K(0,t-s) \cap S^K(0, s)) }_p \\
 &\le&   2 \norm{    m( \widetilde     S^K(0,t-s) \cap S^K(0, \infty)) }_p \\
 &\le& c'_p \, \sqrt{t-s},
\end{eqnarray*}
by Lemma \ref{L:m2b}. Then  (\ref{inc1})  is proven.

For $t\le u \le t+1$, we have \begin{eqnarray*} \norm{ \sup_{t\le
u\le t+1} \big| m(S^K(0, u)) -   m(S^K(0, t)) \big|  }_p  &=&
\norm{\arch{ m(S^K(t,u)
\backslash S^K(0, t))} }_p \\
    & \le & 2 \norm{  m(S^K(t,t+1)
\backslash S^K(0, t))  }_p \\
    & = & 2 \norm{  m(\widetilde S^K(0,1)
\backslash S^K(0, t))  }_p \\
    & \le & 2 \norm{  m(S^K(0,1)) }_p \\
    &\le& c_p,
\end{eqnarray*}

\no which  completes the proof.  \hfill$\Box$

\section{Increments of intersection local times}

We quote  a Tanaka-Rosen-Yor type formula  for the intersection
local time (see \cite{RY91}, \cite{Y85a}, \cite{Y85}):

\begin{lemma}\label{L:alphaat} {\bf (\cite{Y85}, (3.c))} For any $t\ge1$,
we have
$$ \alpha
(0, {\cal A}_t) = - {1\over 2 \pi} \, \int_1^t d B_u \int_0^{u-1}
d s { B_u - B_s \over |B_u - B_s|^3} - {1\over 2 \pi} \,
\int_1^{t-1} d s \left( {1\over |B_t - B_s|} -  {1\over |B_{1+s} -
B_s|} \right).$$
\end{lemma}

We shall estimate the $L^p$ norm of several quantities related to
the intersection local times $\alpha(0, {\cal A}_t)$.  Let us
collect some preliminary results in the following lemma:

\begin{lemma}\label{L:45} Let $t\ge 27$ and
$p\ge2$. There exists  some constant $c_p>0$,
\begin{eqnarray} \norm{ \int_1^t d s {    B_s \over |
B_s|^3} }_p &\le&  c_p \, \sqrt{\log t}\, \label{L45a}  \\
    \norm{ \int_1^t d s {    B_s \over |
B_s|^3} }_2 &= &   \, \sqrt{2 \, \log t}+ O(1), \label{L45b}  \\
    \norm{ \int_0^{t-1}
{ ds \over |B_t - B_s|}}_p &\le & c_p \sqrt t, \label{L5a}
    \\ \norm{ \arch{\int_0^t {
ds \over |B_{1+s} - B_s|}}}_p &\le & c_p \sqrt t. \label{L5b}
   \end{eqnarray}
\end{lemma}

{\no\bf Proof:}  We shall constantly use the
Burkholder-Davis-Gundy inequality for martingales: For any $p>0$,
there exist some constants $c_p>0$ and $c'_p >0$ such that for
any continuous martingale $(M_t)$, we have
\begin{equation}\label{bdg} c_p\, \left( \norm{ \langle M\rangle_t
}_{p/2}\right)^{1/2} \le \norm{ \sup_{0\le s\le t} |M_s|}_p \le
c'_p\,\left( \norm{ \langle M\rangle_t }_{p/2}\right)^{1/2} .
\end{equation}

We define $R_t \df |B_t|, t\ge0$ which is the radial part of the
three dimensional Brownian motion $B$. It is well-known
(\cite{RY99}, Chap. 12) that $R$ is a three-dimensional Bessel
process and admits the following canonical decomposition:
$$ R_t = \gamma_t + {1\over2} \int_0^t { d s \over R_s}, \qquad
t\ge0,$$ for some one-dimensional Brownian motion $\gamma$.

Applying Itô's formula to ${x\over |x|}$ gives that
\begin{eqnarray} \int_1^{t} d s { B_s \over |  B_s|^3}  &= & -
\left( { B_t \over |B_t|} - { B_1\over |B_1|}\right) + \int_1^t {d
B_s \over |B_s|}  - \int_1^t {B_s \over |B_s|}  { B_s \circ d B_s
\over |B_s|^2} \nonumber
    \\ &\df &  - \left( { B_t \over |B_t|} - {
B_1\over |B_1|}\right) + M_t, \label{L5f} \qquad u>1.
    \end{eqnarray}

Using the Burkholder-Davis-Gundy inequality  (\ref{bdg})  for the
martingale $M =(M^{(i)}_u),\, {1\le i\le3}$ implies  that
\begin{eqnarray*} \norm{ M_t }_p &\le & c_p\,  \sum_{i=1}^3 \norm{ \langle
M^{(i)}\rangle_t }^{1/2}_{p/2} \\
    &= & c_p\,  \sum_{i=1}^3 \norm{ \int_1^t d s 
    \left( {1\over |B_s|^2} - {(B^{(i)}_s)^2\over 
    |B_s|^4}\right)}^{1/2}_{p/2}
    \\ &\le& 6 c_p \norm{ \int_1^t {d s\over R^2_s}}^{1/2}_{p/2}
    \\ &\le & c'_p \sqrt{ \log  t},
        \end{eqnarray*}

\no where the last inequality is a consequence of the fact that
\begin{equation} \label{L5a1} {1\over \log t } \int_1^t 
{ ds \over R^2_s} \, \to\, 1,
\qquad \mbox{a.s.  and in  $L^p$}. \end{equation} 
Thus (\ref{L45a}) is proven. Observe that  $$ \e |M_t|^2= \sum_{i=1}^3
\e \langle M^{(i)}\rangle_t = \sum_{i=1}^3 \e \int_1^t ds \left(
{1\over |B_s|^2} - {(B^{(i)}_s)^2\over |B_s|^4}\right) = 2 \log t,
\quad t>1,$$ which together with (\ref{L5f}) yields (\ref{L45b}).
 It follows from Brownian time-reversal that 
\begin{eqnarray*}  \e
\left( \int_0^t { ds \over |B_t - B_s|} \right)^p & = &  \e \left(
\int_0^t { ds \over R_s} \right)^p
    \\& = & 2^p \e\left( R_t -  \gamma_t\right)^p
    \\ &\le & 4^p \left( \e (R_t)^p + \e |\gamma_t|^p  \right)
    \\&=& c'_p t^{p/2},
\end{eqnarray*}

\no by the scaling property. This yields (\ref{L5a}). To obtain
(\ref{L5b}), we write $$ \int_0^t { ds \over |B_{1+s} - B_s|} =
\sum_{j=0}^{[t/2]} \int_{2j}^{2j+1} { ds \over |B_{1+s} - B_s|} +
\sum_{j=0}^{[t/2]} \int_{2j+1}^{2j+2} { ds \over |B_{1+s} - B_s|}
- \int_t^{2[t/2]+2}  { ds \over |B_{1+s} - B_s|}.$$

In each of two sums on $j$, we have independent and identically
distributed variables with common law that of $\int_0^1 { d s
\over |B_{1+s} - B_s|} \in L^p$; hence by Rosenthal's inequality
(\ref{ros}), the two sums centered by their expectations have a
$L^p$ norm bounded by $O(\sqrt{t})$. The $L^p$-norm of the third
term itself is  less than $ \norm{ \int_t^{t+2} { d s \over
|B_{1+s} - B_s|}}_p = \norm{ \int_0^2 { d s \over |B_{1+s} -
B_s|}}_p   $  is  bounded on $t$. Then (\ref{L5b}) is proven and
completes the proof of Lemma \ref{L:45}. \hfill$\Box$

\begin{lemma}\label{L7} For any $\epsilon>0$, almost surely,
\begin{eqnarray} \arch{ \int_1^{t-1} { d s \over |B_{1+s} - B_s|}
  }&=& o\left( \sqrt{t} \, (\log t)^\delta\right),
\qquad {\rm a.s.}, \,\, t\to\infty,  \label{L7a} \\
  \int_1^{t-1} { d  s \over |B_t - B_s|}  
  &=& o\left( \sqrt{t} \, (\log t)^\delta\right),
\qquad {\rm a.s.}, \,\, t\to\infty, \label{L7b}
\end{eqnarray}
\end{lemma}

{\no\bf Proof:} The proof is based on an application of
Propositions \ref{P:inc1} and \ref{P:inc2}. We shall prove that
for any $t>s>2$, \begin{equation} \norm{ \int_0^{t-1}  du {1\over
|B_t - B_u|} - \int_0^{s-1}  du {1\over |B_s -
    B_u|}}_p  \le   c_p  \sqrt{ (t-s) \log t} . \label{L8b}
\end{equation}

\no Therefore,  by virtue of  (\ref{L5a}) and (\ref{L8b}),
Proposition \ref{P:inc2} implies (\ref{L7a}).  The (\ref{L7b}) can
be proven similarly.

To show (\ref{L8b}), we have 
\begin{eqnarray*} && \norm{ \int_0^{t-1}  du {1\over |B_t
- B_u|} - \int_0^{s-1}  du {1\over |B_s - B_u|}}_p
    \\&\le  & \norm{ \int_{s-1}^{t-1}  du
{1\over |B_t - B_u|}  }_p + \norm{ \int_0^{s-1}  du \left({1\over |B_t
- B_u|} - {1\over |B_s - B_u|}\right)}_p
    \\&\le & c_p \sqrt{t-s}+  \norm{ \int_s^t   d B_r \int_0^{s-1}
    d u { B_r - B_u\over |B_r - B_u|^3}}_p,
\end{eqnarray*}

\no by using respectively (\ref{L5a}) to obtain the first term and
the Itô's formula for the second term. Now, by using (\ref{bdg}),
the $L^p$ norm in the above inequality is smaller than
$$   c_p \norm{\int_s^t d r \left| \int_0^{s-1}
    d u { B_r - B_u\over |B_r - B_u|^3}\right|^2 }^{1/2}_{p/2} \le c_p
\left(\int_s^t  dr \norm{\int_1^r { du \over R^2_u}}_p \right)^{1/2}\le
c'_p
    \sqrt{ (t-s) \log t} .$$
\hfill $\Box$

We can now prove our main estimate in this section:

\begin{lemma}\label{L:L8} For $1\le s < t$ and $t>2$,
\begin{eqnarray}  \norm{
\int_1^t du \left| \int_0^{u-1} d s { B_u - B_s \over |B_u
-B_s|^3} \right|^2 - 2 t \log t }_p &\le & c_p \, t \sqrt{ \log t
\, \log\log t} , \label{L5c}
    \\ \norm{ \int_s^t du \left| 
    \int_0^{u-1}  d v { B_u -B_v\over |B_u - B_v|^3} \right|^2 }_p
    &\le & c_p (t-s)  \log  t, \label{L8a}
\end{eqnarray}
\end{lemma}

{\no\bf Proof:}   Notice that by Brownian time-reversal and
(\ref{L45a}),
\begin{equation} \label{L5e} \norm{ \int_0^{u-1} d s { B_u - B_s
\over |B_u - B_s|^3} }_p =\norm{ \int_1^{u} d s {    B_s \over |
B_s|^3} }_p \le c_p \, \sqrt{\log u}\, \qquad u\ge2,
\end{equation}

\no   which yields  (\ref{L8a}). In view of (\ref{L45b}),
(\ref{L5c}) is equivalent to prove that $$ \norm{ \int_1^t du
\arch{ \left| \int_0^{u-1} d s { B_u - B_s \over |B_u -B_s|^3}
\right|^2}  }_p  \le c_p t \sqrt{ \log t \log \log t}.$$

Define $$ Y_u(a,b) = \int_a^b d s { B_u - B_s \over |B_u -B_s|^3},
\qquad 0<a<b<u.$$

Then by scaling property, $\norm{ \int_1^t du \arch{ \left|
\int_0^{u-1} d s { B_u - B_s \over |B_u -B_s|^3} \right|^2}  }_p =
t\, \norm{ \int_{1\over t}^1  d u \arch{ \big| Y_u(0, u-{1\over
t})\big|^2} }_p$. Then it suffices to show that  for $0< \epsilon <
1/27$, \begin{equation} \label{L5g} \norm{ \int_\epsilon^1 d u
\arch{ \big| Y_u(0, u-\epsilon)\big|^2} }_p \le c_p  \sqrt{ \log
{1\over \epsilon} \log \log {1\over \epsilon}}.
\end{equation}

To this end, we remark that for any $0<a<b<u$, \begin{eqnarray*}
\norm{  Y_u(a,b)  }_p &=& \norm{  \int_{u-b}^{u-a} d s {B_s \over
|B_s |^3} }_p\\
    &\le & 2 + \norm{ M_{u-a} - M_{u-b}}_p \\
    & \le & 2 + c_p \norm{ \int_{u-b}^{u-a} {ds \over
    R^2_s}}_{p/2}^{1/2}  \\
    &\le & 2+ c_p \sqrt{ \log { u-a\over u-b}},
\end{eqnarray*}

\no by means of (\ref{L5a1}).  Choose $\delta= {1\over n} \sim
(\log {1\over \epsilon})^{-\theta}$ for some constant $\theta>2$.
We have $$  \int_\epsilon^1 d u \arch{ \big| Y_u(0,
u-\epsilon)\big|^2} =  \int_\epsilon^\delta d u \arch{ \big|
Y_u(0, u-\epsilon)\big|^2} + \int_\delta^1 du \arch{ \left| Y_u(0,
u-\epsilon) \right|^2}.$$

Observe that 
\begin{eqnarray} \norm{\int_\epsilon^\delta d u
\arch{ \big| Y_u(0, u-\epsilon)\big|^2}}_p &\le &
2\norm{\int_\epsilon^\delta d u   \big| Y_u(0, u-\epsilon)\big|^2
}_p  \nonumber\\
    &\le & 2 \int_\epsilon^\delta du \norm{\big| Y_u(0, u-\epsilon)\big|^2
}_p  \nonumber \\
    &\le &  c_p \int_\epsilon^\delta d u (1+ \log (u/\epsilon)) \nonumber
    \\&\le & c_p \delta \log (1/\epsilon)  = c_p (\log
    (1/\epsilon))^{1-\theta}. \label{L5h}
\end{eqnarray}

We claim that \begin{equation}\label{L5i}
 \norm{ \int_\delta^1 du \arch{ \left| Y_u(0,  u-\epsilon) \right|^2} - \int_\delta^1 du
\arch{ \left|  Y_u(u-\delta, u-\epsilon) \right|^2} }_p \le c_p
\sqrt{ \log {1\over \epsilon} \log \log {1\over \epsilon} }.
\end{equation}

In fact, by Cauchy-Schwarz' inequality, $$ \norm{ \int_0^1 d s
f(u) g(u) }_p \le \left( \norm{ \int_0^1 du
|f(u)|^2}_p\right)^{1/2} \left( \norm{ \int_0^1 du
|g(u)|^2}_p\right)^{1/2}.$$

Then by writing $ Y_u(0, u- \epsilon)= Y_u(0, u-\delta) +
Y_u(u-\delta, u-\epsilon)$, we have \begin{eqnarray*}
 && \norm{ \int_\delta^1 du 
\arch{ \left| Y_u(0,  u-\epsilon) \right|^2} - \int_\delta^1 du
\arch{ \left|  Y_u(u-\delta, u-\epsilon) \right|^2} }_p
    \\ &\le & 2  \norm{ \int_\delta^1 du   \left| Y_u(0,  u-\delta)
    \right|^2}_p + 2 \norm{ \int_\delta^1 du   \left| Y_u(0,  u-\delta)
    \right|^2}_p^{1/2} \, \norm{ \int_\delta^1 du   
    \left| Y_u(  u-\delta, u-\epsilon)
    \right|^2}_p^{1/2}.
\end{eqnarray*}

But, $ \norm{ \int_\delta^1 du   \left| Y_u(0,  u-\delta)
    \right|^2}_p  \le c_p \int_\delta^1 du (1+ \log (1/\delta))
    \le c_p \log 1/\delta \le c'p \log \log (1/\epsilon)$, and
    similarly $\norm{ \int_\delta^1 du   \left| Y_u(  u-\delta, u-\epsilon)
    \right|^2}_p \le c_p \int_\delta^1 du (1+ \log
    \delta/\epsilon) \le c'_p \log (1/\epsilon)$. Hence
    (\ref{L5i}) is obtained.

    To bound the $L^p$ norm of $\int_\delta^1 du
\arch{ \left|  Y_u(u-\delta, u-\epsilon) \right|^2}$, we cut the
interval $[\delta, 1]$ into $n=  1/\delta $ parts: $$
\int_\delta^1 du \arch{ \left|  Y_u(u-\delta, u-\epsilon)
\right|^2} = \sum_{j \le n-1: j \mbox { odd}} \int_{j \delta
}^{(j+1) \delta} + \sum_{j \le n-1: j \mbox { even}} \int_{j \delta
}^{(j+1) \delta}.$$

Remark that the sum on odd $j\le n-1$ is a sum of iid variables
whose common law is that of $$ \int_\delta^{2\delta} \arch{ \left|
Y_u(u-\delta, u-\epsilon) \right|^2} \law  \delta\, \int_1^2  du
\arch{ \left| \int_{u-1}^{u - \epsilon/\delta}  d s { B_u - B_s
\over |B_u -B_s|^3} \right|^2} .$$

\no Applying (\ref{ros}), we obtain that \begin{eqnarray*}
\norm{\sum_{j \le n-1: j \mbox { odd}} \int_{j \delta }^{(j+1)
\delta} \left|  Y_u(u-\delta, u-\epsilon) \right|^2  }_p &\le &
c_p\, n^{1/2}\,  \delta\, \norm{\int_1^2  du \arch{ \left|
\int_{u-1}^{u - \epsilon/\delta}  d s { B_u - B_s \over |B_u
-B_s|^3} \right|^2} }_p \\
    &\le & c'_p\, n^{1/2}\,  \delta\, \int_1^2  du (1+ \log
    \delta/\epsilon) \\
    &\le&  c'_p \delta^{1/2}  \log (1/\epsilon) \\
    &=& c'_p (\log 1/\epsilon)^{1-\theta/2}.
\end{eqnarray*}
The same holds for the sum on even $j$. Since $\theta>2$, we get
$\norm{ \int_\delta^1 du \arch{ \left|  Y_u(u-\delta, u-\epsilon)
\right|^2}}_p = O(1)$ which in view  of (\ref{L5h}) and (\ref{L5i})
implies (\ref{L5g}), as desired. \hfill$\Box$

\section{Proof of Theorem \ref{T1}}

Denote by $$ N_t =- {1\over 2 \pi} \, \int_1^t d B_u \int_0^{u-1}
d s { B_u - B_s \over |B_u - B_s|^3}, \qquad t>1.$$

We claim that  
\begin{equation}\label{pra}  X_t \df \arch{
m(S^K(0,t)) } - {\cal C}(K)^2 N_t =   o\left( \sqrt{t} \, (\log
t)^\delta\right), \qquad {\rm a.s.}, \,\, t\to\infty.
\end{equation}

In fact, by Proposition \ref{P:lp}, (\ref{L5a}) and (\ref{L5b}),
we have  \begin{equation}\label{prb}  \norm{ X_t  }_p \le c_p
\sqrt {t}. \end{equation}

On the other hand,   applying (\ref{bdg}) to $N_t$,  we deduce
from Lemma \ref{L:L8} and Corollary \ref{C:inc} that
\begin{eqnarray*} \norm{ X_ t - X_s}_p &\le & c_p \sqrt{(t-s ) \log
t} , \qquad t>s> 27, \\
    \norm{\sup_{t\le u\le t+1} | X_u - X_t| }_p &\le & c_p
    \sqrt{\log t}, \qquad t>27,
\end{eqnarray*}

\no which together with (\ref{prb}) allows us to apply Proposition
\ref{P:inc1}, and shows (\ref{pra}).

We can also apply Proposition \ref{P:inc2} in view of Lemma
\ref{L:L8}, and we obtain that
\begin{eqnarray}
   \int_1^t du \left| \int_0^{u-1} d s { B_u - B_s \over |B_u
-B_s|^3} \right|^2 - 2 t \log t &=& o\left(  t  (\log
t)^{{1\over2}+ \delta}\right), \qquad {\rm a.s.}, \,\, t\to\infty.
\label{pt3}
\end{eqnarray}

By Dubins-Schwarz's representation, there exists some one-dimensional 
Brownian motion $ \beta$ such that  
\begin{eqnarray}
 N_t  &=&  {1\over \pi \sqrt 2} \beta \left( {1\over 2} \int_1^t d u
\left| \int_0^{u-1} d s { B_u - B_s \over |B_u - B_s|^3}
\right|^2\right)  \nonumber \\
    &=&  {1\over \pi \sqrt 2} \beta (    t\log t) +  o\left(  t^{1/2}
(\log t)^{{1\over4}+ \delta}\right), \qquad {\rm a.s.}, \,\,
    t\to\infty, \label{pt4}
\end{eqnarray}

\no by using (\ref{pt3}) and the Brownian increments (cf.
\cite{CR81}, Theorem 1.2.1): For a non-decreasing function $0<
a_t\le t $ such that $t/a_t\uparrow$, we have $$
\limsup_{t\to\infty} {1\over \sqrt{ 2 a_t ( \log (t/a_t)+ \log\log
t)}}\,  \sup_{0\le s \le t- a_t} \sup_{0\le v\le a_t} \, |
\widetilde \beta (s+v) - \widetilde \beta (s)| =1, \qquad \as
$$

Theorem \ref{T1} follows  by assembling (\ref{pra})  and
(\ref{pt4}). \hfill $\Box$

\bigskip
{\noindent\bf Acknowledgements.}  Cooperation between the authors
was supported by the joint French-Hun\-garian Intergovernmental
Grant "Balaton" no. F-39/00.

\end{document}